\documentclass[11pt]{amsart}

\usepackage{amsmath}
\usepackage{amssymb}
\usepackage{amsthm}
\usepackage{mathtools}
\usepackage{tikz}

\usepackage{placeins}

\usepackage{float}

\usepackage{comment}

\theoremstyle{definition}
\newtheorem{thm}{Theorem}[section]

\newtheorem{lem}[thm]{Lemma}

\newtheorem*{ack*}{Acknowledgements}

\subjclass[2010]{05C65(primary)}

\newcommand{\HH}{\mathcal{H}}
\newcommand{\EE}{\mathbb{E}}
\newcommand{\PP}{\mathbb{P}}

\begin{document}

\title{Bollob\'as-Scott Conjecture $3$-uniform}
\author{Hunter Spink, Marius Tiba}

\title{Judiciously $3$-partitioning $3$-uniform hypergraphs}

\begin{abstract}
Bollob\'as, Reed and Thomason proved every $3$-uniform hypergraph $\HH$ with $m$ edges has a vertex-partition $V(\HH)=V_1 \sqcup V_2 \sqcup V_3$ such that each part meets at least $\frac{1}{3}(1-\frac{1}{e})m$ edges, later improved to $0.6m$ by Halsegrave and improved asymptotically to $0.65m+o(m)$ by Ma and Yu. We improve this asymptotic bound to $\frac{19}{27}m+o(m)$, which is best possible up to the error term, resolving a special case of a conjecture of Bollob\'as and Scott.
\end{abstract}

\maketitle

\section{Introduction}
Judicious partitioning problems seek to partition the vertices of a hypergraph $\HH$ such that various quantities are simultaneously maximized. The first such problem was treated by Bollob\'as and Scott \cite{Scott93}, where they proved that one can partition the vertices of a graph with $m$ edges into $k$ parts such that each part contains at most $\frac{m}{k}+o(m)$ edges. They later proved \cite{Scott} that the vertices of a $3$-uniform hypergraph can be partitioned into $k$ parts each of which contains at most $\frac{m}{k^2}+o(m)$ edges. Recently Hou, Wu, Zeng and Zhu \cite{Hou} claimed to have shown the result for $4$-uniform hypergraphs with $\frac{m}{k^3}+o(m)$ edges, but their key technical lemma is incorrect.

The judicious partitioning problem we consider involves partitioning the vertices of an $r$-uniform hypergraph with $m$ edges into $k$ parts so that the minimum number of edges touched by any part is maximized. Bollob\'as and Scott \cite{Scott3} conjectured that this maximum is $(1-(1-\frac{1}{k})^r)m+o(m)$, which if true is optimal up to the error term by considering the complete $r$-uniform hypergraph. Special cases of this problem have garnered considerable interest. Bollob\'as, Reed and Thomason \cite{BRT} proved every $3$-uniform hypergraph $\HH$ has a vertex-partition $V(\HH)=V_1 \sqcup V_2 \sqcup V_3$ such that each part meets at least $\frac{1}{3}(1-\frac{1}{e})m$ edges, which was claimed to have been improved to $\frac{5m-1}{9}$ by Bollob\'as and Scott \cite{Scott3} (although their proof apparently contains a subtle error, see Halsegrave \cite{Hal}), improved to $0.6m$ by Halsegrave \cite{Hal} and improved asymptotically to $0.65m+o(m)$ by Ma and Yu \cite{Ma}.

In this article, we solve the conjecture of Bollob\'as and Scott in this special case.
\begin{thm}
Every $3$-uniform hypergraph $\HH$ has a vertex-partition $V(\HH)=V_1 \sqcup V_2 \sqcup V_3$ such that each part meets at least $\frac{19}{27}m-O(m^{6/7})$ edges.
\end{thm}

Our proof follows a strategy of Bollob\'as and Scott \cite{Scott} and Ma and Yu \cite{Ma}. The key technical lemma in \cite{Ma} cannot be optimized sufficiently to attain $\frac{19}{27}$, in fact cannot be pushed past $0.7<\frac{19}{27}=0.\overline{703}$. We make an additional observation which imposes an additional inequality in this technical lemma.

We now outline the strategy. It would be nice if for a given $3$-uniform hypergraph, the uniformly random $3$-partition of the vertices worked. Each part would touch in expectation $\frac{19}{27}m+o(m)$ edges, but the presence of high degree vertices prevents these three quantities from simultaneously concentrating around the mean. For example, if there is a pair of vertices which belong to every edge, then a part which contains neither vertex will touch with high probability only $\frac{1}{3}m+o(m)$ edges. To avoid this situation, we partition the vertices according to a two stage process.

We first partition ``high'' degree vertices $V^{high}\subseteq V(\HH)$ into three parts $V_1,V_2,V_3$, and then take the remaining ``low'' degree vertices $V^{low}=V(\HH)\setminus V^{high}$ and assign each of these vertices independently to $V_i$ with some probability $p_i$ with $p_1+p_2+p_3=1$, so that the expected number of edges touched in each part is at least $\frac{19}{27}m+o(m)$. Because we exclude the high degree vertices from the random process, we have a tight concentration around the mean in each part by the Azuma-Hoeffding inequality (or more specifically the version of the inequality due to McDiarmid \cite{McDiarmid}), and we conclude.

The threshold for ``high'' degree vertices and the application of the Azuma-Hoeffding inequality are essentially automatic (if we only desire an $o(m)$ error, we only have to ensure that the threshold for ``high'' degree is $o(m)$). The challenge is to appropriately partition the ``high'' degree vertices so that the $p_i$ can be chosen to make each part touch the correct number of edges in expectation.

As in \cite{Ma}, after choosing the high degree threshold, we partition $V^{high}$ as follows. Denote by $e(\HH)^i$ for $i \in \{0,1,2,3\}$ to be the edges in $e(\HH)$ which touch $i$ vertices in $V^{high}$ and $3-i$ vertices in $V^{low}$. The edges in $e(\HH)^2$ induce a multi-graph $G$ with vertex set $V^{high}$ obtained by replacing each edge $e$ with $e \cap V^{high}$. Denote by $x_i$ the number of edges of $G$ within $V^{high}_i$, and $b_{ij}$ the number of edges of $G$ with endpoints in $V^{high}_i$ and $V^{high}_j$. We choose our partition $V^{high}=V^{high}_1 \sqcup V^{high}_2 \sqcup V^{high}_3$ so that $b_{23}+b_{13}+b_{12}$ is maximal (equivalently $x_1+x_2+x_3$ is minimal).  The only facts used in \cite{Ma} about this partition were that $b_{ij} \ge \max(2x_i, 2x_j)$. These inequalites follow by noting that the multigraph between $V^{high}_i$ and $V^{high}_j$ maximizes the number of edges between the two parts for fixed $V^{high}_i \sqcup V^{high}_j$.

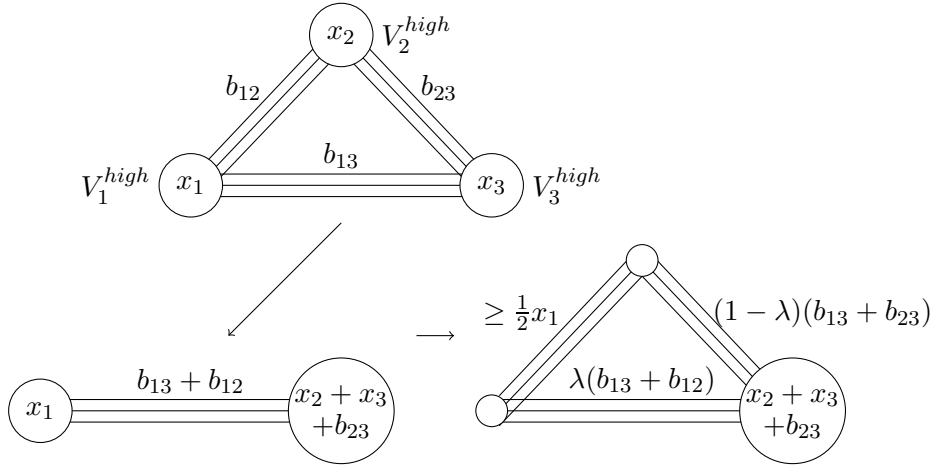
\begin{center}
\begin{tikzpicture}
\draw (-1.9,0.15)--(1.9,0.15);
\draw (-1.9,0)--(1.9,0);
\draw (-1.9,-0.15)--(1.9,-0.15);
\draw (-1.9,0.2)--(-0.2,2);
\draw (-1.9,0)--(0,2);
\draw (-1.9,-0.2)--(0.2,2);
\draw (1.9,0.2)--(0.2,2);
\draw (1.9,0)--(0,2);
\draw (1.9,-0.2)--(-0.2,2);
\draw (-1.3,1.3) node{$b_{12}$};
\draw (1.3,1.3) node{$b_{23}$};
\draw (0,0.4) node{$b_{13}$};
\filldraw [draw=black, fill=white] (-2,0) circle (12pt);
\filldraw [draw=black, fill=white] (2,0) circle (12pt);
\filldraw [draw=black, fill=white] (0,2) circle (12pt);
\draw (-3,0) node{$V^{high}_1$};
\draw (-2,0) node{$x_1$};
\draw (1,2) node{$V^{high}_2$};
\draw (0,2) node{$x_2$};
\draw (2,0) node{$x_3$};
\draw (3,0) node{$V^{high}_3$};
\draw[->] (0,-0.5)--(-1.5,-2);
\draw (-1.9-2,0.15-3)--(1.9-2,0.15-3);
\draw (-1.9-2,0-3)--(1.9-2,0-3);
\draw (-1.9-2,-0.15-3)--(1.9-2,-0.15-3);
\filldraw [draw=black, fill=white] (-4,-3) circle (12pt);
\filldraw [draw=black, fill=white] (-0,-3) circle (20pt);
\draw (-4,-3) node{$x_1$};
\draw (-0,-2.8) node{$x_2+x_3$};
\draw (-0,-3.2) node{$+b_{23}$};
\draw (-2,-2.6) node{$b_{13}+b_{12}$};
\draw[->] (1,-2)--(1.5,-2);
\def\t{4}
\def\s{3}
\draw (-1.9+\t,0.15-\s)--(1.9+\t,0.15-\s);
\draw (-1.9+\t,0-\s)--(1.9+\t,0-\s);
\draw (-1.9+\t,-0.15-\s)--(1.9+\t,-0.15-\s);
\draw (-1.9+\t,0.2-\s)--(-0.2+\t,2-\s);
\draw (-1.9+\t,0-\s)--(0+\t,2-\s);
\draw (-1.9+\t,-0.2-\s)--(0.2+\t,2-\s);
\draw (1.9+\t,0.2-\s)--(0.2+\t,2-\s);
\draw (1.9+\t,0-\s)--(0+\t,2-\s);
\draw (1.9+\t,-0.2-\s)--(-0.2+\t,2-\s);
\draw (-1.6+\t,1.3-\s) node{$\ge \frac{1}{2}x_1$};
\draw (2.4+\t,1.3-\s) node{$(1-\lambda)(b_{13}+b_{23})$};
\draw (0+\t,0.4-\s) node{$\lambda (b_{13}+b_{12})$};
\filldraw [draw=black, fill=white] (-2+\t,0-\s) circle (6pt);
\filldraw [draw=black, fill=white] (2+\t,0-\s) circle (20pt);
\filldraw [draw=black, fill=white] (0+\t,2-\s) circle (6pt);
\draw (2+\t,0.2-\s) node{$x_2+x_3$};
\draw (2+\t,-0.2-\s) node{$+b_{23}$};
\end{tikzpicture}
\end{center}

However, these numerical constraints fail to rule out the possibility that $G$ is (close to) a configuration which is complete in $V^{high}_k$, and complete bipartite between $V^{high}_k$ and the other two parts, which prevents the method from proving a bound exceeding $0.7m$. We identify the additional inequality that if $G$ maximizes $b_{23}+b_{13}+b_{12}$, then $b_{ij} \ge \frac{1}{2}x_k$, where $k$ is the index $\ne i,j$. This arises from considering the effect of moving all vertices from $V^{high}_i$ to $V^{high}_j$, decreasing $b_{23}+b_{13}+b_{12}$ by $b_{ij}$, and dividing $V^{high}_k$ between $V^{high}_k$ and $V^{high}_i$ in such a way that $b_{ik}$ is maximized, increasing $b_{23}+b_{13}+b_{12}$ by at least $\frac{1}{2}x_k$.

The core of this paper is a rather technical lemma. This technical lemma is used to show that if we subdivide $V^{high}$ so that $b_{23}+b_{13}+b_{12}$ is maximal, then we can find probabilities $p_1,p_2,p_3$ as described earlier. Given an instance of $\HH$, finding appropriate values of $p_i$ amount to simply solving cubic equations. Hence the laborious method used to prove the inequality is not reflected in the elegance of the final algorithm. We would be very interested in a less opaque derivation of the technical lemma.

The technical lemma is a system of inequalities with $13$ variables and $11$ degrees of freedom. We ultimately reduce this system to systems of inequalities with $6$ variables and $2$ degrees of freedom. For the systems which are not tight, we finely subdivide our parameter space and use a computer to check the system on each box. We prove all remaining systems of inequalities in the body of the paper.
\section{The technical lemma, its contrapositive, and our strategy}
The following technical lemma is the crux of our proof.
\begin{lem}
For any non-negative variables $$x_1+x_2+x_3+a_1+a_2+a_3+b_{12}+b_{13}+b_{23}+c=1$$ 
with $$b_{ij} \ge \max(2x_i,2x_j,\frac{1}{2}x_k),$$ there exists $q_1,q_2,q_3 \in [0,1]$ with $q_1+q_2+q_3=2$ such that the following inequalties are satisfied.
\begin{align*}
q_1(b_{23}+x_2+x_3)+q_1^2(a_2+a_3)+q_1^3c \le \frac{8}{27}\\
q_2(b_{13}+x_1+x_3)+q_2^2(a_1+a_3)+q_2^3c \le \frac{8}{27}\\
q_3(b_{12}+x_1+x_2)+q_3^2(a_1+a_2)+q_3^3c \le \frac{8}{27}.
\end{align*}
\end{lem}
This lemma improves the lemma in \cite{Ma}, which does not impose the inequality $b_{ij} \ge \frac{1}{2}x_k$, and has the constant $\frac{8}{27}$ relaxed $0.35$. The following three sections are devoted to the proof of this lemma.

As the left hand sides of the inequalities in the system are non-decreasing in the $q_i$ and are zero when $q_i=0$, by symmetry the contrapositive of the lemma is equivalent to the following.

There does not exist $q_1,q_2,q_3 \in [0,1]$ with $q_1+q_2+q_3=2$ and non-negative variables $$x_1+x_2+x_3+b_{23}+b_{13}+b_{12}+a_1+a_2+a_3+c=1$$
with $$b_{ij} \ge \max(2x_i,2x_j,\frac{1}{2}x_k)$$ such that one of the following two systems holds.

\textbf{System1}
\begin{align*}
q_1(b_{23}+x_2+x_3)+q_1^2(a_2+a_3)+q_1^3c > \frac{8}{27}\\
q_2(b_{13}+x_1+x_3)+q_2^2(a_1+a_3)+q_2^3c > \frac{8}{27}\\
q_3(b_{12}+x_1+x_2)+q_3^2(a_1+a_2)+q_3^3c > \frac{8}{27}.
\end{align*}

\textbf{System2} $q_1=1$ (so $q_2+q_3=1$) and
\begin{align*}
q_2(b_{13}+x_1+x_3)+q_2^2(a_1+a_3)+q_2^3c > \frac{8}{27}\\
q_3(b_{12}+x_1+x_2)+q_3^2(a_1+a_2)+q_3^3c > \frac{8}{27}.
\end{align*}

We will arrive at a contradiction assuming one of the two systems holds.

We may replace $(a_1,a_2,a_3,c) \mapsto (a_1+(1-q_1)c,a_2+(1-q_2)c,a_3+(1-q_3)c,0)$ and the corresponding system still holds, so we may assume from now on that $c=0$. By working with the contrapositive we may assume $q_1,q_2,q_3$ are fixed and hence we are able to eliminate $c$ in this fashion, which was not exploited previously, but drastically reduces the algebraic complexity of all future systems of inequalities.

Denoting $L=x_1+x_2+x_3+b_{23}+b_{13}+b_{12}+a_1+a_2+a_3+c$ and $L_i=q_i(b_{jk}+x_j+x_k)+q_i^2(x_j+x_k)+q_i^3c$, our strategy for \textbf{System1} is as follows. Keeping the $q_i$ fixed, we repeatedly linearly perturb the remaining variables in such a way that $L$ stays constant and $L_1,L_2,L_3$ do not decrease. We then perturb in this direction until one of the constraints attains equality. We then merge the associated variables (or set the corresponding variable to $0$ if a variable hits zero), and repeat.

By linear algebra, it is clear we can do this (e.g. by choosing any perturbation fixing $L,L_1,L_2$ and choosing the sign so that $L_3$ does not decrease), until the number of variables besides $q_1,q_2,q_3$ is $3$. However, we do some manual reductions first to reduce the final number of cases and simplify the exposition.

We handle \textbf{System2} in Section 3. \textbf{System1} will be dealt with in Sections 4, 5, and 6, completing the proof of Lemma 2.1.

\section{System2}
Recall we may assume $c=0$.
We may now replace
\begin{align*}
(x_1,x_2,x_3,a_1,a_2,a_3)\mapsto (x_1,0&,0,0,0,0)\\
(b_{23},b_{13},b_{12})\mapsto(\frac{1}{2}x_1,&b_{13}+x_3+q_2(a_1+a_2+a_3),\\&b_{12}+x_2+q_3(a_1+a_2+a_3)+b_{23}-\frac{1}{2}x_1)
\end{align*}
so we may assume that $x_2=x_3=a_1=a_2=a_3=0$ and $b_{23}=\frac{1}{2}x_1$. The system is now $b_{12},b_{13} \ge 2x_1 \ge 0$, $$b_{12}+b_{13}+\frac{3}{2}x_1=1$$ and
\begin{align*}
q_2(b_{13}+x_1) > \frac{8}{27}\\
q_3(b_{12}+x_1)> \frac{8}{27}.
\end{align*}
As $b_{12},b_{13}\ge 2x_1$, we have $x_1 \le \frac{2}{11}$. Hence,
\begin{align*}
1&>\frac{8}{27}(\frac{1}{b_{13}+x_1}+\frac{1}{b_{12}+x_1})\\
&\ge \frac{8}{27}\frac{(1+1)^2}{b_{13}+b_{23}+2x_1}\\
&=\frac{32}{27}\frac{1}{1+\frac{1}{2}x_1} \ge \frac{32}{27}\frac{1}{1+\frac{1}{11}}=\frac{88}{81},
\end{align*}
a contradiction.

\section{System1 assuming $x_1\ge 4x_2$ and $x_2\ge x_3$}
In this section, we prove the result assuming $x_1 \ge 4x_2$ and $x_2 \ge x_3$. The inequalities governing the $b_{ij}$ are thus $b_{12},b_{13}\ge 2x_1$ and $b_{23}\ge \frac{1}{2}x_1$. Our plan is to vary the variables, with the promise that we never break any of the constraints of \textbf{System1}. Recall we may assume $c=0$. Then we have $$x_1+x_2+x_3+b_{23}+b_{13}+b_{12}+a_1+a_2+a_3=1$$ and
\begin{align*}
q_1(b_{23}+x_2+x_3)+q_1^2(a_2+a_3) > \frac{8}{27}\\
q_2(b_{13}+x_1+x_3)+q_2^2(a_1+a_3) > \frac{8}{27}\\
q_3(b_{12}+x_1+x_2)+q_3^2(a_1+a_2) > \frac{8}{27}.
\end{align*}

\textbf{Cases 1a and 1b} First, suppose that $x_2=x_3$. Call this common value $x_{23}$. Then we have $$x_1+2x_{23}+b_{23}+b_{13}+b_{12}+a_1+a_2+a_3=1$$ and
\begin{align*}
q_1(b_{23}+2x_{23})+q_1^2(a_2+a_3) > \frac{8}{27}\\
q_2(b_{13}+x_1+x_{23})+q_2^2(a_1+a_3) > \frac{8}{27}\\
q_3(b_{12}+x_1+x_{23})+q_3^2(a_1+a_2) > \frac{8}{27}.
\end{align*}
Now, replace $(x_{23},b_{23})\mapsto (x_{23},b_{23})+\epsilon(1,-2)$ starting with $\epsilon=0$, and increase $\epsilon$ until either $x_1=4x_{23}$ or $b_{23}=\frac{1}{2}x_1$.

\textbf{Case1a}
Suppose first that $x_1=4x_{23}$, so that $b_{12},b_{13}\ge 8x_{23}$, $b_{23} \ge 2x_{23}$, and $x_{23}\ge 0$. Then we have $$6x_{23}+b_{23}+b_{13}+b_{12}+a_1+a_2+a_3=1$$ and
\begin{align*}
q_1(b_{23}+2x_{23})+q_1^2(a_2+a_3) > \frac{8}{27}\\
q_2(b_{13}+5x_{23})+q_2^2(a_1+a_3) > \frac{8}{27}\\
q_3(b_{12}+5x_{23})+q_3^2(a_1+a_2) > \frac{8}{27}.
\end{align*}
Now, replace $(x_{23},b_{23},b_{13},b_{12})\mapsto (x_{23},b_{23},b_{13},b_{12})+\epsilon(1,2,-4,-4)$ starting with $\epsilon=0$, and increase $\epsilon$ until either $b_{12}=8x_{23}$ or $b_{13}=8x_{23}$. Without loss of generality assume $b_{12}=8x_{23}$. Then we have

\textbf{System1a} $b_{13} \ge 8x_{23}$ and $b_{23} \ge 2x_{23}$ and $x_{23} \ge 0$ and $$14x_{23}+b_{23}+b_{13}+a_1+a_2+a_3=1$$ and
\begin{align*}
q_1(b_{23}+2x_{23})+q_1^2(a_2+a_3)> \frac{8}{27}\\
q_2(b_{13}+5x_{23})+q_2^2(a_1+a_3)> \frac{8}{27}\\
q_3(13x_{23})+q_3^2(a_1+a_2)> \frac{8}{27}.
\end{align*}

\textbf{Case1b}
Suppose now instead that $b_{23}=\frac{1}{2}x_1$. Then we have $$\frac{3}{2}x_1+2x_{23}+b_{12}+b_{13}+a_1+a_2+a_3=1$$ and
\begin{align*}
q_1(\frac{1}{2}x_1+2x_{23})+q_1^2(a_2+a_3)>\frac{8}{27}\\
q_2(b_{13}+x_1+x_{23})+q_2^2(a_1+a_3)>\frac{8}{27}\\
q_3(b_{12}+x_1+x_{23})+q_3^2(a_1+a_2)>\frac{8}{27}.
\end{align*}
Now, replace $(x_1,b_{12},b_{13})\mapsto (x_1,b_{12},b_{13})+\epsilon(1,-\frac{3}{4},-\frac{3}{4})$ starting with $\epsilon=0$, and increase $\epsilon$ until either $b_{12}=2x_1$ or $b_{13}=2x_1$. Without loss of generality assume $b_{12}=2x_1$. Then we have

\textbf{System1b} $b_{13}\ge 2x_1 \ge 8x_{23}\ge 0$ and $$\frac{7}{2}x_1+2x_{23}+b_{13}+a_1+a_2+a_3=1$$ and
\begin{align*}
q_1(\frac{1}{2}x_1+2x_{23})+q_1^2(a_2+a_3)>\frac{8}{27}\\
q_2(b_{13}+x_1+x_{23})+q_2^2(a_1+a_3)>\frac{8}{27}\\
q_3(3x_1+x_{23})+q_3(a_1+a_2)>\frac{8}{27}.
\end{align*}

\textbf{Cases 1c and 1d}
Now, we suppose that $x_2 \ge x_3$ and equality is not attained. We will do similar replacements of variables, but if we ever have $x_2=x_3$ then we can reduce to either \textbf{System1a} or \textbf{System1b}, so we may suppose this never happens.

Replace $(x_3,b_{23})\mapsto (x_3,b_{23})+\epsilon(1,-1)$ starting with $\epsilon=0$, and increase $\epsilon$ until $b_{23}=\frac{1}{2}x_1$. Now, replace $(x_3,b_{13})\mapsto (x_3,b_{13})+\epsilon(1,-1)$ starting with $\epsilon=0$, and increase $\epsilon$ until $b_{13}=2x_1$. Then we have $$\frac{7}{2}x_1+x_2+x_3+b_{12}+a_1+a_2+a_3=1$$ and
\begin{align*}
q_1(\frac{1}{2}x_1+x_2+x_3)+q_1^2(a_2+a_3)>\frac{8}{27}\\
q_2(3x_1+x_3)+q_2^2(a_1+a_3)>\frac{8}{27}\\
q_3(b_{12}+x_1+x_2)+q_3^2(a_1+a_2)>\frac{8}{27}.
\end{align*}
Replace $(x_2,b_{12})\mapsto (x_2,b_{12})+\epsilon(1,-1)$ starting with $\epsilon=0$, and increase $\epsilon$ until either $b_{12}=2x_1$ or $x_1=4x_2$.

\textbf{Case1c}
Suppose first that $b_{12}=2x_1$. Then we have

\textbf{System1c} $x_1\ge x_2 \ge x_3 \ge 0$ and $$\frac{11}{2}x_1+x_2+x_3+a_1+a_2+a_3=1$$ and
\begin{align*}
q_1(\frac{1}{2}x_1+x_2+x_3)+q_1^2(a_2+a_3)>\frac{8}{27}\\
q_2(3x_1+x_3)+q_2^2(a_1+a_3)>\frac{8}{27}\\
q_3(3x_1+x_2)+q_3^2(a_1+a_2)>\frac{8}{27}.
\end{align*}

\textbf{Case1d}
Finally, suppose that $x_1=4x_2$. Then we have

\textbf{System1d} $b_{12} \ge 8x_2 \ge 8x_3 \ge 0$ and $$15x_2+x_3+b_{12}+a_1+a_2+a_3=1$$ and
\begin{align*}
q_1(3x_2+x_3)+q_1^2(a_2+a_3)>\frac{8}{27}\\
q_2(12x_2+x_3)+q_2^2(a_1+a_3)>\frac{8}{27}\\
q_3(b_{12}+5x_2)+q_3^2(a_1+a_3)>\frac{8}{27}.
\end{align*}

\subsection{Solving System1a,b,c,d}
For each of these systems, we carry out the procedure described at the end of Section 2. This process terminates when we have $3$ variables apart from $q_1,q_2,q_3$, and for each system there are $\binom{6}{3}=20$ resulting cases. For each of \textbf{System1a,b,c,d} one of the 20 resulting cases will be when $a_1=a_2=a_3=0$, and one of the resulting cases will be when the remaining 3 variables are equal to zero. 

In the former case, when $a_1=a_2=a_3=0$, we have
\begin{align*}
2&>\frac{8}{27}(\frac{1}{b_{23}+x_2+x_3}+\frac{1}{b_{13}+x_1+x_3}+\frac{1}{b_{12}+x_1+x_2}) \\
&\ge \frac{8}{27}\frac{(1+1+1)^2}{b_{23}+x_2+x_3+b_{13}+x_1+x_3+b_{12}+x_1+x_2}\\
&=\frac{8}{3(1+x_1+x_2+x_3)}\ge\frac{8}{3(1+\frac{1}{3})}=2,
\end{align*}
a contradiction.

In the latter case, we have $a_1+a_2+a_3=1$ and $$2> \sqrt{\frac{8}{27}}(\frac{1}{\sqrt{a_2+a_3}}+\frac{1}{\sqrt{a_1+a_3}}+\frac{1}{\sqrt{a_1+a_2}}).$$ The right hand side is at least $2$ by Jensen's inequality applied to $f(x)=\frac{1}{\sqrt{x}}$, a contradiction.

Aside from these two cases for each system, we handle all other cases by computer as follows. Let $\tilde{q_i}$ be the minimum of $1$ and the solution to the quadratic equation $L_i=\frac{8}{27}$ (recall that $c=0$). Then we may equivalently check that $\tilde{q_1}+\tilde{q_2}+\tilde{q_3}\ge 2$ for any choice of the remaining variables. As it turns out, the infimum of $\tilde{q_1}+\tilde{q_2}+\tilde{q_3}$ is strictly larger than $2$ in all remaining cases, which we can check by subdividing the space of the remaining variables very finely (boxes of dimension $0.001\times 0.001 \times 0.001$ suffices in all cases), and crudely bounding in each box the smallest possible value of $\tilde{q_1}+\tilde{q_2}+\tilde{q_3}$. We report the results of these computations in the Appendix.

\section{System1 assuming $4x_2 \ge x_1 \ge x_2 \ge x_3$}

We assume from now on that $4x_2 \ge x_1 \ge x_2 \ge x_3$, and so the inequalities with the $b_{ij}$ become $b_{12},b_{13} \ge 2x_1$, $b_{23} \ge 2x_2$. The symmetry of \textbf{System1}, along with the previous two sections, completes the verification of the technical lemma. As in the previous section, we will use the method described at the end of Section 2. However, if we ever attain equality in $4x_2\ge x_1$, then we may stop and use Section 4, so we may assume from now on that we never attain equality in $4x_2 \ge x_1$. Recall we may assume $c=0$. Then we have $$x_1+x_2+x_3+b_{23}+b_{13}+b_{12}+a_1+a_2+a_3=1$$ and
\begin{align*}
q_1(b_{23}+x_2+x_3)+q_1^2(a_2+a_3) > \frac{8}{27}\\
q_2(b_{13}+x_1+x_3)+q_2^2(a_1+a_3) > \frac{8}{27}\\
q_3(b_{12}+x_1+x_2)+q_3^2(a_1+a_2) > \frac{8}{27}.
\end{align*}

Replace $(x_1,x_2,x_3,b_{23},b_{13},b_{12}) \mapsto (x_1,x_2,x_3,b_{23},b_{13},b_{12})+\epsilon(-1,-1,3,-2,-2,3)$ starting with $\epsilon=0$, and increase $\epsilon$ until $x_2=x_3$. Call this common value $x_{23}$, so that $b_{12},b_{13} \ge 2x_1$, $b_{23} \ge 2x_{23}$ and $4x_{23} \ge x_1 \ge x_{23}$. Then we have $$x_1+2x_{23}+b_{23}+b_{13}+b_{12}+a_1+a_2+a_3=1$$ and
\begin{align*}
q_1(b_{23}+2x_{23})+q_1^2(a_2+a_3) > \frac{8}{27}\\
q_2(b_{13}+x_1+x_{23})+q_2^2(a_1+a_3) > \frac{8}{27}\\
q_3(b_{12}+x_1+x_{23})+q_3^2(a_1+a_2) > \frac{8}{27}.
\end{align*}
We now replace $(x_1,b_{13})\mapsto (x_1,b_{13})+\epsilon(1,-1)$ starting with $\epsilon=0$, and increase $\epsilon$ until either $b_{13}=2x_1$ or $b_{12}=2x_1$. Without loss of generality, assume that $b_{13}=2x_1$. Finally, replace $(x_{23},b_{23})\mapsto (x_{23},b_{23})+\epsilon(1,-2)$ starting with $\epsilon=0$, and increase $\epsilon$ until either $x_{23}=x_1$ or $b_{23}=2x_{23}$.

In the former case, let $x$ be the common value of $x_{23}$ and $x_1$. Set $A=4x$, $B=b_{23}+2x$, $C=b_{12}+2x$. Then (assuming $C \ge B$ without loss of generality) we have

\textbf{System1e} $C \ge B \ge A \ge 0$ and $$\frac{1}{4}A+B+C+a_1+a_2+a_3=1$$ and
\begin{align*}
q_1B+q_1^2(a_2+a_3) > \frac{8}{27}\\
q_2A+q_2^2(a_1+a_3) > \frac{8}{27}\\
q_3C+q_3^2(a_1+a_2) > \frac{8}{27}.
\end{align*}
In the latter case, set $A=4x_{23}$, $B=3x_1+x_{23}$, and $C=b_{12}+x_1+x_{23}$. Then we have

\textbf{System1f} $C \ge B \ge A \ge 0$ and $B \le \frac{13}{4}A$ and $$\frac{7}{12}A+\frac{2}{3}B+C+a_1+a_2+a_3=1$$ and
\begin{align*}
q_1A+q_1^2(a_2+a_3) > \frac{8}{27}\\
q_2B+q_2^2(a_1+a_3) > \frac{8}{27}\\
q_3C+q_3^2(a_1+a_2) > \frac{8}{27}.
\end{align*}

For each of these systems, we carry out the procedure described at the end of Section 2. This process terminates when we have $3$ variables apart from $q_1,q_2,q_3$. Note in \textbf{System1f} that if we have $B=\frac{13}{4}A$ then $x_1=4x_{23}$ so this is handled by Section 4 and we can avoid analyzing all of these cases. Hence we are left with $20$ cases for each system. The case $A=B=C=0$ is handled as in Section 4.1.

Every case when $A=0$ other than $A=B=C=0$ we check by computer as before. Note that \textbf{System1e} is identical to \textbf{System1f} when $A=B$ as $1+\frac{1}{4}=\frac{2}{3}+\frac{7}{12}$, so for these cases we only need to check one of the two systems. We report the results of these computations in the Appendix.

\subsection{Remaining cases of Systems1e,f}
The only cases remaining are the 20 cases, 10 from \textbf{System1e} and 10 \textbf{System1f}, when three of $A=B,B=C, a_1=0,a_2=0,a_3=0$ occur. Note that $\frac{7}{12}A+\frac{2}{3}B+C+a_1+a_2+a_3 \le \frac{1}{4}A+B+C+a_1+a_2+a_3$ when $A \le B \le C$, so it suffices to check the 10 cases as above in the following system.

\textbf{System1f'} $C \ge B \ge A \ge 0$ and $$\frac{7}{12}A+\frac{2}{3}B+C+a_1+a_2+a_3=1$$ and
\begin{align*}
q_1A+q_1^2(a_2+a_3) > \frac{8}{27}\\
q_2B+q_2^2(a_1+a_3) > \frac{8}{27}\\
q_3C+q_3^2(a_1+a_2) > \frac{8}{27}.
\end{align*}

It is vitally important to note that there are solutions to \textbf{System1f'}, but we claim there are no solutions in the 10 cases mentioned above.

Suppose there is a solution in one of the cases above, we will derive a contradiction. We have
\begin{align*}
2=&q_1+q_2+q_3\\
>&\frac{16}{27}\frac{1}{A+\sqrt{A^2+\frac{32}{27}(a_2+a_3)}}+\frac{1}{B+\sqrt{B^2+\frac{32}{27}(a_1+a_3)}}+ \frac{1}{C+\sqrt{C^2+\frac{32}{27}(a_1+a_2)}}\\
\ge& \frac{16}{3}\frac{1}{A+B+C+\sqrt{A^2+\frac{32}{27}(a_2+a_3)}+\sqrt{B^2+\frac{32}{27}(a_1+a_3)}+\sqrt{C^2+\frac{32}{27}(a_1+a_2)}}
\end{align*}
and therefore
\begin{align*}
A+B+C+\sqrt{A^2+\frac{32}{27}(a_2+a_3)}+\sqrt{B^2+\frac{32}{27}(a_1+a_3)}+\sqrt{C^2+\frac{32}{27}(a_1+a_2)} > \frac{8}{3}.
\end{align*}
We will in fact show that $$A+B+C+\sqrt{A^2+\frac{32}{27}(a_2+a_3)}+\sqrt{B^2+\frac{32}{27}(a_1+a_3)}+\sqrt{C^2+\frac{32}{27}(a_1+a_2)} \le \frac{8}{3}$$ in each of the 10 cases, yielding a contradiction.

Suppose first we are in one of the three cases where $A=B=C$. Then $\frac{9}{4}A+a_1+a_2+a_3=1$, and by Jensen's inequality applied to $\sqrt{x}$, it suffices to show that
$$3A+\sqrt{3}\sqrt{3A^2+\frac{64}{27}(a_1+a_2+a_3)} \le \frac{8}{3}.$$ As $A \le \frac{4}{9}$, we have $\frac{8}{3} \ge 3A$. As $a_1+a_2+a_3=1-\frac{9}{4}A$, this inequality is $$3A+|\frac{8}{3}-3A| \le \frac{8}{3},$$
and as $\frac{8}{3} \ge 3A$, both sides are equal.

If $a_1=a_2=a_3=0$, then we have $2(A+B+C) \le \frac{8}{3}(\frac{7}{12}A+\frac{2}{3}B+C)=\frac{8}{3}$.

In the remaining cases, we have one of $A=B$ and $B=C$. We have dealt with all cases where the maximum occurs on the boundary except when $A=0$, so we can assume the maximum is attained in the interior of the domain or $A=0$.

When $A=B$, we have by symmetry two distinct cases, when $a_1=a_2=0$ and when $a_1=a_3=0$.

When $a_1=a_2=0$ we have $0 \le A \le C$, $\frac{5}{4}A+C+a_3=1$ and we want to show
$$2A+2C+2\sqrt{A^2+\frac{32}{27}a_3} \le \frac{8}{3}.$$

By Lagrange multipliers, if the triple $(A,C,a_3)$ attains a maximum on the interior of the domain, we have
\begin{align*}
\frac{4}{5}(2+\frac{2A}{\sqrt{A^2+\frac{32}{27}a_3}})&=2\\
\frac{32}{27}\frac{1}{\sqrt{A^2+\frac{32}{27}a_3}}&=2.
\end{align*}
But then $\sqrt{A^2+\frac{32}{27}a_3}=\frac{16}{27}$, so $A=\frac{4}{27}$, $a_3=\frac{5}{18}$, and so $C=\frac{29}{54}$. These values make the left hand side $\frac{23}{9}<\frac{8}{3}$. If the maximum occurs when $A=0$, then the inequality is equivalent to $(\frac{4}{3}-C)^2 \ge \frac{32}{27}(1-C)$, which is true.\\

When $a_1=a_3=0$ we have $0 \le A \le C$, $\frac{5}{4}A+C+a_2=1$, and we want to show
$$3A+C+\sqrt{A^2+\frac{32}{27}a_2}+\sqrt{C^2+\frac{32}{27}a_2} \le \frac{8}{3}.$$
If the maximum is attained in the interior of the domain, by Lagrange multipliers
\begin{align*}
3+\frac{A}{\sqrt{A^2+\frac{32}{27}a_2}}=\frac{5}{4}+\frac{5}{4}\frac{C}{\sqrt{C^2+\frac{32}{27}a_2}},
\end{align*}
which cannot happen as $\frac{C}{\sqrt{C^2+\frac{32}{37}a_2}} \le 1$. If the maximum occurs when $A=0$, then the left hand side is at most $C+\sqrt{2C^2+\frac{128}{27}a_2}$, and the inequality then becomes $(\frac{8}{3}-C)^2\ge 2C^2+\frac{128}{27}(1-C)$, which is true.\\

Now, when $B=C$, we have by symmetry two distinct cases, when $a_2=a_3=0$ and when $a_1=a_3=0$.

When $a_2=a_3=0$, we have $0 \le A \le C$, $\frac{7}{12}A+\frac{5}{3}C+a_1=1$, and we want to show
$$2A+2C+2\sqrt{C^2+\frac{32}{27}a_1} \le \frac{8}{3}.$$
If the maximum is attained in the interior of the domain, by Lagrange multipliers
\begin{align*}
-\frac{120}{21}+2+\frac{2C}{\sqrt{C^2+\frac{32}{27}a_1}}=0
\end{align*}
which cannot happen as $\frac{2C}{\sqrt{C^2+\frac{32}{27}a_1}} \le 2$. If the maximum occurs when $A=0$, then the inequality is equivalent to $(\frac{4}{3}-C)^2 \ge C^2+\frac{32}{27}(1-\frac{5}{3}C)$, which is true as $C \le \frac{3}{5}$.\\

Finally, when $a_1=a_3=0$, we have $0 \le A \le C$, $\frac{7}{12}A+\frac{5}{3}C+a_2=1$, and we want to show
$$A+3C+\sqrt{A^2+\frac{32}{27}a_2}+\sqrt{C^2+\frac{32}{27}a_2} \le \frac{8}{3}.$$
By Jensen's inequality applied to $\sqrt{x}$, it suffices to show that
\begin{align*}
A+3C+\sqrt{2A^2+2C^2+\frac{128}{27}a_2} \le \frac{8}{3}.
\end{align*}
As $A+3C \le \frac{8}{3}(\frac{7}{12}A+\frac{5}{3}C) \le \frac{8}{3}$, we can move $A+3C$ to the right hand side, square the equation, and then replace $a_2$ with $1-\frac{7}{12}A-\frac{5}{3}C$ to get
$$(A+3C-\frac{8}{3})^2-(2A^2+2C^2+\frac{128}{27}(1-\frac{7}{12}A-\frac{5}{3}C)) \ge 0.$$
If the minimum is attained in the interior of the domain, then we have the partial derivatives with respect to $A$ and $C$ vanishing. This yields
\begin{align*}
2(A+3C-\frac{8}{3})-4A+\frac{224}{81}&=0\\
6(A+3C-\frac{8}{3})-4C+\frac{640}{81}&=0,
\end{align*}
so $A=\frac{16}{81},C=\frac{40}{81}$. Evaluating at this point yields $\frac{256}{2187} \ge 0.$ The boundary case $A=C$ is identical to the case $A=B=C$ when $a_1=a_3=0$. The boundary case $A=0$ is equivalent to $(\frac{8}{3}-3C)^2 \ge 2C^2+\frac{128}{27}(1-\frac{5}{3}c)$ for $0 \le C \le \frac{3}{5}$, and the boundary case $a_2=0$ is equivalent to $(\frac{8}{3}-3C-(\frac{12}{7}-\frac{20}{7}C))^2\ge(2C^2+2(\frac{12}{7}-\frac{20}{7}C)^2)$ for $\frac{4}{9} \le C \le \frac{3}{5}$, which are both true.

\section{Proof of Theorem 1.1}
In this section, we prove Theorem 1.1 using the Lemma 2.1.

We can assume every vertex has positive degree. Order the vertices $v_1,\ldots,v_n$ such that $d(v_1) \ge d(v_2) \ge \ldots \ge d(v_n)$. Let $V^{high}$ consist of the $t=m^{\alpha}$ vertices of highest degree with some fixed $0<\alpha<\frac{1}{3}$. Then
$$m^\alpha d(v_{t+1}) < \sum_{v \in \HH} d(v)=3m,$$ so $$d(v_{t+1}) \le 3m^{1-\alpha}.$$
Hence in particular, $$\sum_{i=t+1}^n d(v_i)^2 < 3m^{1-\alpha} \sum_{i=t+1}^n d(v_i) \le 9m^{2-\alpha}.$$
Recall $x_1,x_2,x_3,b_{12},b_{13},b_{23}$ were defined in the introduction. Further define $a_i$ for $i \in \{1,2,3\}$ to be the number of edges in $e(\HH)^1$ such that the high degree vertex lies in $V^{high}_i$, and let $c$ be the number of edges in $e(\HH)^0$.

Apply Lemma 2.1 to all of our variables scaled by $\frac{1}{m-e(V^{high})}$ to get probabilities $p_i=1-q_i$. Place the low degree vertices independently into $V_i$ with probability $p_i$. Then $\EE(d(V_i)) \ge \frac{19}{27}(m-e(V^{high}))$.

We recall that the Azuma-Hoeffding inequality \cite{McDiarmid} asserts that if $X_1,\ldots,X_n$ are independent random variables with values in $\{1,\ldots,k\}$, and $f:\{1,\ldots,k\}^n \to \mathbb{N}$ such that $|f(Y)-f(Y')| \le d_i$ if $Y$ and $Y'$ differ only on the $i$'th coordinate, then for $t\ge 0$ we have $\PP(f(X_1,\ldots,X_n) \le \EE(f(X_1,\ldots,X_n))-t) \le e^{-\frac{2t^2}{\sum d_i^2}}$. 

By the Azuma-Hoeffding inequality with $k=3$, we have $$\PP(d(V_i)<\EE(d(V_i))-z) \le e^{\frac{-z^2}{\sum_{i=t+1}^n d(v_i)^2}} < e^{-\frac{2z^2}{9m^{2-\alpha}}}.$$
Taking $z=\sqrt{(9/2) \ln 3}m^{1-\frac{\alpha}{2}}$ yields $\PP(d(V_i)<\EE(d(V_i))-z) < \frac{1}{3}$. Hence there is a partition which has $d(V_i)>\frac{19}{27}(m-e(V^{high}))-z$. $e(V^{high}) = O(m^{3\alpha})$ and $z=O(m^{1-\frac{\alpha}{2}})$, so taking $\alpha=\frac{2}{7},$ we have that the result with $\frac{19}{27}m-O(m^{6/7})$.

\begin{ack*}
The first author would like to thank Harvard university for travel funding, and the second author would like to thank the Trinity Hall Research Studentship. We would both like to thank our advisor B\'ela Bollob\'as for all of his support.
\end{ack*}

\appendix
\section{Results of Computation}
Tables 1--6 show the results of computation. We denote by $\epsilon$ the side length of the cubes used in the interval method (when applicable). We use $\ast$ to indicate that the minimum value was proved in the body of the paper. We record a proven lower bound in the ``Lower Bound'' column.  Finally, we note that when $B=A$ in \textbf{System1f}, the case is identical to the corresponding case in \textbf{System1e}.
\FloatBarrier
\begin{table}
\caption{System1a}
\centering
\resizebox{290pt}{!}{
\begin{tabular}{|c|c|c|c|c|c|c|c|}
\hline
$x_{23}=0$ & $b_{13}=8x_{23}$ & $b_{23}=2x_{23}$ & $a_1=0$ & $a_2=0$ & $a_3=0$ &$\epsilon$&Lower Bound\\
\hline
\checkmark & \checkmark & \checkmark & & & & &$2\ast$\\ \hline
\checkmark & \checkmark & &\checkmark & & &0.002 &2.078\\ \hline
\checkmark & \checkmark & & &\checkmark & &0.002 &2.077\\ \hline
\checkmark & \checkmark & & & &\checkmark &0.002 &2.077\\ \hline
\checkmark &  &\checkmark &\checkmark & & &0.002 &2.077\\ \hline
\checkmark &  &\checkmark & &\checkmark & &0.002 &2.078\\ \hline
\checkmark &  &\checkmark & & &\checkmark &0.002 &2.077\\ \hline
\checkmark &  & &\checkmark &\checkmark & &0.002 &2.085\\ \hline
\checkmark &  & &\checkmark & &\checkmark &0.002 &2.086\\ \hline
\checkmark &  & & &\checkmark &\checkmark &0.002 &2.086\\ \hline
&\checkmark  &\checkmark &\checkmark & & &0.001 &2.005\\ \hline
&\checkmark  &\checkmark & &\checkmark & &0.002 &2.033\\ \hline
&\checkmark  &\checkmark & & &\checkmark &0.002 &2.033\\ \hline
&\checkmark  & &\checkmark &\checkmark & &0.002 &2.057\\ \hline
&\checkmark  & &\checkmark & &\checkmark &0.002 &2.057\\ \hline
&\checkmark  & & &\checkmark &\checkmark &0.002 &2.043\\ \hline
& &\checkmark&\checkmark  &\checkmark &&0.002 &2.045\\ \hline
& &\checkmark&\checkmark  & &\checkmark&0.002 &2.044\\ \hline
& &\checkmark&  &\checkmark &\checkmark&0.002 &2.041\\ \hline
& &&\checkmark  &\checkmark &\checkmark&0.002&$2\ast$ or 2.046\\ \hline
\end{tabular}
}
\end{table}

\begin{table}
\caption{System1b}
\centering
\resizebox{290pt}{!}{
\begin{tabular}{|c|c|c|c|c|c|c|c|}
\hline
$x_{23}=0$ & $b_{13}=2x_1$ & $x_{1}=4x_{23}$ &  $a_1=0$ & $a_2=0$ & $a_3=0$ &$\epsilon$&Lower Bound\\
\hline
\checkmark & \checkmark & \checkmark & & & & &$2\ast$\\ \hline
\checkmark & \checkmark & &\checkmark & & &0.002 &2.042\\ \hline
\checkmark & \checkmark & & &\checkmark & &0.002 &2.069\\ \hline
\checkmark & \checkmark & & & &\checkmark &0.002&2.069\\ \hline
\checkmark &  &\checkmark &\checkmark & & &0.002 &2.077\\ \hline
\checkmark &  &\checkmark & &\checkmark & &0.002 &2.078\\ \hline
\checkmark &  &\checkmark & & &\checkmark &0.002 &2.077\\ \hline
\checkmark &  & &\checkmark &\checkmark & &0.002 &2.072\\ \hline
\checkmark &  & &\checkmark & &\checkmark &0.002 &2.072\\ \hline
\checkmark &  & & &\checkmark &\checkmark &0.002 &2.070\\ \hline
&\checkmark  &\checkmark &\checkmark & & &0.001 &2.005\\ \hline
&\checkmark  &\checkmark & &\checkmark & &0.002 &2.033\\ \hline
&\checkmark  &\checkmark & & &\checkmark &0.002 &2.033\\ \hline
&\checkmark  & &\checkmark &\checkmark & &0.002 &2.026\\ \hline
&\checkmark  & &\checkmark & &\checkmark &0.002 &2.026\\ \hline
&\checkmark  & & &\checkmark &\checkmark &0.002 &2.024\\ \hline
& &\checkmark&\checkmark  &\checkmark &&0.002 &2.045\\ \hline
& &\checkmark&\checkmark  & &\checkmark&0.002 &2.044\\ \hline
& &\checkmark&  &\checkmark &\checkmark&0.002 &2.041\\ \hline
& &&\checkmark  &\checkmark &\checkmark&0.002&$2\ast$ or 2.025\\ \hline
\end{tabular}
}
\end{table}

\begin{table}
\caption{System1c}
\centering
\resizebox{300pt}{!}{
\begin{tabular}{|c|c|c|c|c|c|c|c|}
\hline
$x_3=0$ & $x_2=x_3$ & $x_2=x_1$ & $a_1=0$ & $a_2=0$ & $a_3=0$ &$\epsilon$&Lower Bound\\
\hline
\checkmark & \checkmark & \checkmark & & & & &$2\ast$\\ \hline
\checkmark & \checkmark & &\checkmark & & &0.002 &2.042\\ \hline
\checkmark & \checkmark & & &\checkmark & &0.002 &2.069\\ \hline
\checkmark & \checkmark & & & &\checkmark &0.002 &2.069\\ \hline
\checkmark &  &\checkmark &\checkmark & & &0.001 &2.019\\ \hline
\checkmark &  &\checkmark & &\checkmark & &0.002 &2.036\\ \hline
\checkmark &  &\checkmark & & &\checkmark &0.002 &2.037\\ \hline
\checkmark &  & &\checkmark &\checkmark & &0.002 &2.027\\ \hline
\checkmark &  & &\checkmark & &\checkmark &0.002 &2.028\\ \hline
\checkmark &  & & &\checkmark &\checkmark &0.002 &2.026\\ \hline
&\checkmark  &\checkmark &\checkmark & & &0.001 &2.005\\ \hline
&\checkmark  &\checkmark & &\checkmark & &0.002 &2.033\\ \hline
&\checkmark  &\checkmark & & &\checkmark &0.002 &2.033\\ \hline
&\checkmark  & &\checkmark &\checkmark & &0.002 &2.026\\ \hline
&\checkmark  & &\checkmark & &\checkmark &0.002 &2.026\\ \hline
&\checkmark  & & &\checkmark &\checkmark &0.002 &2.024\\ \hline
& &\checkmark&\checkmark  &\checkmark &&0.001 &2.042\\ \hline
& &\checkmark&\checkmark  & &\checkmark&0.001 &2.042\\ \hline
& &\checkmark&  &\checkmark &\checkmark&0.001 &2.042\\ \hline
& &&\checkmark  &\checkmark &\checkmark&0.001 &$2\ast$ or 2.033\\ \hline
\end{tabular}
}
\end{table}

\begin{table}
\caption{System1d}
\centering
\resizebox{300pt}{!}{
\begin{tabular}{|c|c|c|c|c|c|c|c|}
\hline
$b_{12}=8x_{2}$ & $x_2=x_3$ & $x_3=0$ & $a_1=0$ & $a_2=0$ & $a_3=0$ &$\epsilon$&Lower Bound\\
\hline
\checkmark & \checkmark & \checkmark & & & & &$2\ast$\\ \hline
\checkmark & \checkmark & &\checkmark & & &0.002 &2.077\\ \hline
\checkmark & \checkmark & & &\checkmark & &0.002 &2.077\\ \hline
\checkmark & \checkmark & & & &\checkmark &0.002&2.078\\ \hline
\checkmark &  &\checkmark &\checkmark & & &0.001 &2.019\\ \hline
\checkmark &  &\checkmark & &\checkmark & &0.002 &2.036\\ \hline
\checkmark &  &\checkmark & & &\checkmark &0.002 &2.037\\ \hline
\checkmark &  & &\checkmark &\checkmark & &0.002 &2.047\\ \hline
\checkmark &  & &\checkmark & &\checkmark &0.002 &2.047\\ \hline
\checkmark &  & & &\checkmark &\checkmark &0.002 &2.041\\ \hline
&\checkmark  &\checkmark &\checkmark & & &0.001 &2.005\\ \hline
&\checkmark  &\checkmark & &\checkmark & &0.002 &2.033\\ \hline
&\checkmark  &\checkmark & & &\checkmark &0.002 &2.033\\ \hline
&\checkmark  & &\checkmark &\checkmark & &0.002 &2.044\\ \hline
&\checkmark  & &\checkmark & &\checkmark &0.002 &2.045\\ \hline
&\checkmark  & & &\checkmark &\checkmark &0.002 &2.041\\ \hline
& &\checkmark&\checkmark  &\checkmark &&0.001 &2.042\\ \hline
& &\checkmark&\checkmark  & &\checkmark&0.001 &2.042\\ \hline
& &\checkmark&  &\checkmark &\checkmark&0.001 &2.042\\ \hline
& &&\checkmark  &\checkmark &\checkmark&0.001 &$2\ast$ or 2.042\\ \hline
\end{tabular}
}
\end{table}

\begin{table}
\caption{System1e}
\centering
\resizebox{300pt}{!}{
\begin{tabular}{|c|c|c|c|c|c|c|c|}
\hline
$A=0$ & $B=A$ & $C=B$ & $a_1=0$ & $a_2=0$ & $a_3=0$ &$\epsilon$&Lower Bound\\
\hline
\checkmark & \checkmark & \checkmark & & & & &$2\ast$\\ \hline
\checkmark & \checkmark & &\checkmark & & &0.002 &2.077\\ \hline
\checkmark & \checkmark & & &\checkmark & &0.002 &2.077\\ \hline
\checkmark & \checkmark & & & &\checkmark &0.002 &2.078\\ \hline
\checkmark &  &\checkmark &\checkmark & & &0.002 &2.075\\ \hline
\checkmark &  &\checkmark & &\checkmark & &0.002 &2.076\\ \hline
\checkmark &  &\checkmark & & &\checkmark &0.002 &2.076\\ \hline
\checkmark &  & &\checkmark &\checkmark & &0.002 &2.086\\ \hline
\checkmark &  & &\checkmark & &\checkmark &0.002 &2.085\\ \hline
\checkmark &  & & &\checkmark &\checkmark &0.002 &2.084\\ \hline
&\checkmark  &\checkmark &\checkmark & & & &$2\ast$\\ \hline
&\checkmark  &\checkmark & &\checkmark & & &$2\ast$\\ \hline
&\checkmark  &\checkmark & & &\checkmark & &$2\ast$\\ \hline
&\checkmark  & &\checkmark &\checkmark & & &$2\ast$\\ \hline
&\checkmark  & &\checkmark & &\checkmark & &$2\ast$\\ \hline
&\checkmark  & & &\checkmark &\checkmark & &$2\ast$\\ \hline
& &\checkmark&\checkmark  &\checkmark && &$2\ast$\\ \hline
& &\checkmark&\checkmark  & &\checkmark& &$2\ast$\\ \hline
& &\checkmark&  &\checkmark &\checkmark& &$2\ast$\\ \hline
& &&\checkmark  &\checkmark &\checkmark& &$2\ast$\\ \hline
\end{tabular}
}
\end{table}

\begin{table}
\caption{System1f}
\centering
\resizebox{300pt}{!}{
\begin{tabular}{|c|c|c|c|c|c|c|c|}
\hline
$A=0$ & $B=A$ & $C=B$ & $a_1=0$ & $a_2=0$ & $a_3=0$ &$\epsilon$&Lower Bound\\
\hline
\checkmark & \checkmark & \checkmark & & & & &$2\ast$\\ \hline
\checkmark & \checkmark & &\checkmark & & &0.002 &2.077\\ \hline
\checkmark & \checkmark & & &\checkmark & &0.002 &2.077\\ \hline
\checkmark & \checkmark & & & &\checkmark &0.002 &2.078\\ \hline
\checkmark &  &\checkmark &\checkmark & & &0.002 &2.075\\ \hline
\checkmark &  &\checkmark & &\checkmark & &0.002 &2.076\\ \hline
\checkmark &  &\checkmark & & &\checkmark &0.002 &2.076\\ \hline
\checkmark &  & &\checkmark &\checkmark & &0.002 &2.086\\ \hline
\checkmark &  & &\checkmark & &\checkmark &0.002 &2.085\\ \hline
\checkmark &  & & &\checkmark &\checkmark &0.002 &2.084\\ \hline
&\checkmark  &\checkmark &\checkmark & & & &$2\ast$\\ \hline
&\checkmark  &\checkmark & &\checkmark & & &$2\ast$\\ \hline
&\checkmark  &\checkmark & & &\checkmark & &$2\ast$\\ \hline
&\checkmark  & &\checkmark &\checkmark & & &$2\ast$\\ \hline
&\checkmark  & &\checkmark & &\checkmark & &$2\ast$\\ \hline
&\checkmark  & & &\checkmark &\checkmark & &$2\ast$\\ \hline
& &\checkmark&\checkmark  &\checkmark && &$2\ast$\\ \hline
& &\checkmark&\checkmark  & &\checkmark& &$2\ast$\\ \hline
& &\checkmark&  &\checkmark &\checkmark& &$2\ast$\\ \hline
& &&\checkmark  &\checkmark &\checkmark& &$2\ast$\\ \hline
\end{tabular}
}
\end{table}

\FloatBarrier

\end{document}